\documentclass[11pt]{amsart}
\usepackage{amsmath,amsfonts,amssymb,amsthm,amscd,latexsym}
\usepackage{hyperref}
\usepackage{xcolor}
\usepackage{enumitem}
\usepackage{chngcntr}
\usepackage{pst-all}
\begin{document}
\newtheorem{theorem}{Theorem}[section]
\newtheorem{definition}[theorem]{Definition}
\newtheorem{proposition}[theorem]{Proposition}
\newtheorem{lemma}[theorem]{Lemma}
\newtheorem{remark}[theorem]{Remark}
\newtheorem{corollary}[theorem]{Corollary}
\newtheorem{question}{Question}
\newtheorem{example}{Examples}[section]
\newtheorem{notation}[theorem]{Notation}
\newtheorem{claim}[theorem]{Claim}
\newtheorem{fact}[theorem]{Fact}
\newcommand\cl{\begin{claim}}
\newcommand\ecl{\end{claim}}
\newcommand\rem{\begin{remark}\upshape}
\newcommand\erem{\end{remark}}
\newcommand\ex{\begin{example}\upshape}
\newcommand\eex{\end{example}}
\newcommand\nota{\begin{notation}\upshape}
\newcommand\enota{\end{notation}}
\newcommand\dfn{\begin{definition}\upshape}
\newcommand\edfn{\end{definition}}
\newcommand\cor{\begin{corollary}}
\newcommand\ecor{\end{corollary}}
\newcommand\thm{\begin{theorem}}
\newcommand\ethm{\end{theorem}}
\newcommand\prop{\begin{proposition}}
\newcommand\eprop{\end{proposition}}
\newcommand\lem{\begin{lemma}}
\newcommand\elem{\end{lemma}}
\newcommand\fct{\begin{fact}}
\newcommand\efct{\end{fact}}
\providecommand\qed{\hfill$\quad\Box$}
\newcommand\pr{{Proof:\;}}
\newcommand\prcl{\par\noindent{\em Proof of Claim: }}
\newcommand{\R}{\mathbb R} 
\newcommand{\Q}{\mathbb Q}  
\newcommand{\Z}{\mathbb Z}
\newcommand{\bF}{\mathbb F}
\newcommand{\bN}{\mathbb N} 
\newcommand{\cM}{\mathcal M}
\newcommand{\cB}{\mathcal B}
\newcommand{\cL}{\mathcal L}
\newcommand{\cC}{\mathcal C}
\newcommand{\cA}{\mathcal A}
\newcommand{\cK}{\mathcal K}
\newcommand{\cD}{\mathcal D}
\newcommand{\cR}{\mathcal R}
\newcommand{\cV}{\mathcal V}
\newcommand{\cG}{\mathcal G}
\newcommand{\cP}{\mathcal P}
\newcommand{\bC}{\mathfrak C}
\newcommand{\Pp}{\mathbb P}
\newcommand{\Aa}{\mathbb A}
\def\union{\cup}
\newcommand{\pf}{\noindent{\em Proof. }}
\newcommand{\la}{\langle}
\newcommand{\ra}{\rangle}
\newcommand\cX{{{\mathcal{X}}}}
\newcommand\y{{{\mathbf{y}}}}
\newcommand\x{{{\mathbf{x}}}}
\newcommand\ba{{{\mathbf{a}}}}
\newcommand\bb{{{\mathbf{b}}}}
\newcommand\cc{{{\mathbf{c}}}}
\newcommand\bu{{{\mathbf{u}}}}
\newcommand{\acl}{\mathrm{acl}}
\newcommand{\dcl}{\mathrm{dcl}}
\newcommand{\Cl}{\mathrm{cl}}
\newcommand{\divi}{\mathrm{div}}
\newcommand{\Div}{\mathrm{Div}}

\newcommand{\im}{\mathrm{im}}
\newcommand{\I}{\mathrm{I}}
\newcommand{\pCF}{\mathrm{pCF}}
\newcommand{\RCVF}{\mathrm{RCVF}}
\newcommand{\ACVF}{\mathrm{ACVF}}
\newcommand{\IFTA}{\mathrm{IFTA}}
\newcommand{\RCF}{\mathrm{RCF}}
\newcommand{\RV}{\mathrm{RV}}

\newcommand{\ord}{\text{ord}}
\newcommand{\Ja}{\nabla}

\author{Fran\c coise Point}
\email{point@math.univ-paris-diderot.fr}
\title{Dense pairs of rings}

\date{\today}
\begin{abstract}  
Outside of the framework of geometric theories, we exhibit complete, respectively model-complete theories of rings whose corresponding theory of pairs 
is complete, respectively model-complete, using transfer results proven in the seventies for boolean products of structures. It includes certain boolean products of pairs of dp-minimal fields of characteristic $0$. We also show, as in the case of pairs of fields \cite{CP}, how it fits in the framework of differential rings.
\end{abstract}
\maketitle
\section{Introduction}

\par The study of theories of dense pairs of structures have a long history, starting with a  result of A. Robinson on the completeness and model-completeness of the theory of dense pairs of real-closed fields \cite{R} and most of the further developments took place in the general framework of lovely pairs of geometric theories \cite{BV}, or dense pairs of complete theories with an existential matroid extending the theory of integral domains \cite{F}.

However C. Toffalori  \cite{T} proved a transfer result using the result mentioned above of A. Robinson to theories of boolean products over a space with no isolated points, in the same line of transfer of completeness and or model-completeness of (certain) theories of fields to theories of von Neumann regular rings \cite{LS}, \cite{M}.
\par Here, we generalize Toffalori's result in order to describe theories of dense pairs of boolean products of certain geometric theories of fields, in particular open theories of topological fields \cite{CP} (see section \ref{top}).
We place ourselves in the framework developed by S. Burris and H. Werner which encompasses a number of previous transfer of first-order properties in these products such as completeness and model-completeness \cite{BW}. We use completeness of theories of dense pairs of geometric theories of fields \cite{F}. We also propose a definition of a dense pair of boolean products of geometric theories of fields (see section \ref{top}).
\par Then we also want to show, as in the case of pairs of fields, how it fits in the framework of differential rings \cite[section  4]{CP}, using this time transfer results in boolean products of theories of differential fields. 
\par In \cite{CP}, we showed in some specific setting that certain differential expansions of NIP theories of fields do retain the NIP property. Of course one cannot expect it here in these boolean products of differential fields, but in case the differential fields have a NIP theory, the so-called determining sequence of a formula (in the Feferman-Vaught theorem on products), consists of a formula in the language of boolean algebra and finitely many NIP-formulas (in the theory of factors).
\par Note that recent works address the question, in the case of rings, which constraints on the algebraic structure, combinatorial model-theoretic properties such as NIP or dp-minimality or having finite dp-rank impose. 
For instance a NIP ring has finitely many maximal ideals  \cite[Proposition 2.1]{EH}, dp-minimal integral domain is a local ring and a ring of finite dp-rank is a direct product of finitely many henselian local rings \cite{J}. 
\par One can apply our results to obtain complete, respectively model-complete theories of certain boolean products of dense pairs of dp-minimal fields of characteristic $0$, using former results of W. Johnson \cite{J23} on the algebraic structure of dp-minimal fields (see section \ref{app}). 
\section{Boolean products}
Let $\cC$ be a class of $\cL$-structures.
 We consider subdirect products $\cA:=\prod_{x\in X}^s \cA_{x}$ of elements $\cA_{x}\in \cC$ over some index set $X$, namely $\cL$-substructures of direct products of elements of $\cC$ with the additional property that for any $x\in X$, and any $a_{x}\in A_{x}$, there is $a:=(a(y))_{y\in X}\in \cA$ such that $a(x)=a_{x}$. When all structures $\cA_{x}$ are the same, say $\cD$, we denote the direct product (over $X$) by $\cD^X$. 
\nota Let $\cC$ be a class of $\cL$-structures and let $\cA:=\prod_{x\in X}^s \cA_{x}$ be a subdirect product of elements $\cA_{x}\in \cC$ over some index set $X$. Let $\varphi(x_{1},\ldots,x_{n})$ be an $\cL$-formula and let $\bar f:=(f_{1},\ldots,f_{n})\in \prod_{x\in X}^s \cA_{x}$. Then $[\varphi(\bar f)]:=\{x\in X\colon \cA_{x}\models \varphi(f_{1}(x),\ldots,f_{n}(x))\}$.
\enota
\par Recall that by Stone representation theorem, any boolean algebra is isomorphic to the boolean algebra of continuous functions on a totally disconnected compact Hausdorff space $\cX$ (also called boolean space) with values in the boolean ring $\Z/2\Z$.
\par Let $\cX$ a Boolean space and denote by $\cX^*$ the boolean algebra of clopen subsets of $\cX$. 

Let $\cD\in \cC$ and consider the following subdirect product in $\cD^{X}$: the $\cL$-substructure $\cD[\cX]^*$ whose domain consists of $\{f\in D^{X}\colon f^{-1}(d)$ is a clopen subset of $\cX$ for every $d\in D\}$; it is called the bounded boolean power of $\cD$.
\medskip
 \par Fix $\tilde T$ a theory of Boolean algebras. Burris and Werner considered the following classes  $\Gamma_{\tilde T}^a(\cC)$, $\Gamma_{\tilde T}^e(\cC)$ of subdirect products of elements of $\cC$ \cite[section 1]{BW}.
\dfn \label{bp} The class $\Gamma_{\tilde T}^a(\cC)$ of $\cL$-structures consists of all subdirect products $\cA=\prod^s_{x\in X}\cA_{x}$ of elements of $\cC$, with $X=\cX$ a boolean space and $\cX^*\models \tilde T$, which in addition satisfy the following:
\begin{enumerate}[label={(P\arabic*)}]
\item atomic extension property: if $\varphi$ is an atomic $\cL$-formula, for any $\bar f\in A$, then $[\varphi(\bar f)]$ is a clopen subset of $\cX$,
\item patchwork property: for any $f, g\in A$ and any clopen subset $U$ of $\cX$, there is $h\in A$ such that $U\subseteq[f=h]$ and $X\setminus U\subseteq[g=h]$.
\end{enumerate}
The subclass $\Gamma_{\tilde T}^e(\cC)$ of $\Gamma_{\tilde T}^e(\cC)$ consists of those elements of $\Gamma_{\tilde T}^a(\cC)$ which satisfy in addition:
\begin{enumerate}[label={(P\arabic*)}, resume]
\item elementary extension property: if $\varphi$ is an $\cL$-formula and $\bar f\in A$, then $[\varphi(\bar f)]$ is a clopen subset of $\cX$.
\end{enumerate}
For $\cA\in \Gamma_{\tilde T}^a(\cC)$, we will denote by $\cX(A)$ the underlying boolean space (or by $\cX$ if this is no ambiguity) and sometimes we will omit the subscript $\tilde T$ (when it does not play a role).
\edfn
\par Let $\cD\in \cC$ and $\cX^*\models \tilde T$, then the bounded boolean power $\cD[\cX]^*$ 
belongs to $\Gamma_{\tilde T}^e(\cC)$ \cite[section 2]{BW}. Furthermore, when $\cC$ is the class of models of a complete $\cL$-theory $T$, then for any $\cD\in \cC$, $Th(\Gamma^e_{\tilde T}(\cC))=Th\{\cD[\cX]^*\colon \cX^*\models \tilde T\}$ \cite[Theorem 4.5(c)]{BW}. We will denote $Th(\Gamma^e_{\tilde T}(\cC))$ by $\Gamma^e_{\tilde T}(T)$. In particular $Th(\Gamma^e_{\tilde T}(\cC))$ is complete whenever $\tilde T$ is complete. The main tool in proving this result is a Feferman-Vaught theorem on products, revisited for sheaves of structures by, for instance S. Comer, and stated in \cite{BW} as follows. 
\par To every $\cL$-formula $\varphi$, first one can (effectively) associate a determining sequence, namely a sequence of formulas consisting of a formula $\Phi^*(z_{1},\ldots,z_{\ell})$ in the language of boolean algebras and finitely many $\cL$-formulas $\psi_{1},\ldots, \psi_{\ell}$. Then, this determining sequence allows one to reduce satisfaction in certain products to satisfaction in the factors and in the underlying boolean algebra.
\par The determining sequence is constructed by induction on the complexity of $\varphi$: one puts $\varphi$ in prenex normal form and one describes the procedure first for atomic formulas, then how to proceed with negations and conjunctions and finally with formulas with one existential quantifier.
\fct \cite[Theorem 4.1]{BW}\label{det} Let $\varphi(\bar u)$ be an $\cL$-formula with a determining sequence $$(\Phi^*(z_{1},\ldots,z_{\ell}),\psi_{1}(\bar u),\ldots,\psi_{\ell}(\bar u)).$$
Then for $\cA\in \Gamma^e(\cC)$ and $\bar f\in A$, we have
\[
\cA\models \varphi(\bar f) \leftrightarrow \cX(A)^*\models \Phi^*([\psi_{1}(\bar f)],\ldots,[\psi_{\ell}(\bar f)]).
\]
\efct
\section{Discriminators and existentially closed boolean products}\label{projector}
\par When $\cC$ is a class of $\cL$-structures with a model-complete theory, one can look for conditions implying that an element of $\Gamma^a_{\tilde T}(\cC)$ belongs to $\Gamma^e_{\tilde T}(\cC)$. In that perspective one enriches the language $\cL$ with a discriminator \cite[section 9]{BW}. Since we are interested in classes $\cC$ of $\cL$-structures which expand an abelian group, assuming that the language $\cL$  contain the group language $\{+,-, 0\}$, we instead introduce a projector, namely a binary function symbol $p(u,v)$ defined by $p(a,b)=a$ if $b=0$ and $p(a,b)=0$ otherwise. This binary function will be used as a discriminator (in this particular class of structures), namely a term with the property that $t(u,v,w)=z$ if and only if $(u=v\;\wedge \;w=z)\; \vee\; (u\neq v\; \wedge\; u=z)$ \cite[section 9]{BW}. (So in our setting, the following discriminator formula is $p(u-w,u-v)=u-z$.)
\medskip
\par Denote by $\cC^{p}$ the class of expansions of elements of $\cC$ in the language $\cL_{p}:=\cL\cup\{p(.,.)\}$.

\par The existence of a discriminator formula allows one to axiomatize the existentially closed boolean products of models of a model-complete $\cL$-theory $T$ in $\cL_{p}$ \cite[Theorem 10.7]{BW}. (Note that in Theorem 10.7 in \cite{BW}, one assumes that the language only contains function symbols, but later on (see pages 305--306 in \cite{BW}) the authors give conditions of the language in order to handle the case where $\cL$ contains relation symbols (and get the analog of Theorem 10.7). 
 They require that for each $n$-ary relation $r(\bar x)$, there is a positive existential $\cL$-formula $\varphi_{r}(\bar x)$ such that 
\[
\;\; (\dagger)\;\;\;\;T\models \forall \bar x (\neg r(\bar x)\leftrightarrow \varphi_{r}(\bar x)).
\]

\par In particular one can check that the analog of \cite[Lemma 9.1]{BW} holds in this setting. 
\lem Then in $\cC$ we have:
\begin{enumerate}
\item $\cC^{p}\models (u=0\vee v= 0) \leftrightarrow p(u,v)=u$,
\item $\cC^{p}\models (u=0\wedge v= 0) \leftrightarrow p(u,v)+v=0$,
\item $\cC^{p}\models (u=0\vee v\neq 0) \leftrightarrow p(u,v)=0$.
\end{enumerate}
In particular, in $\cC^{p}$, any open $\cL$-formula not containing relation symbols, is equivalent either to an atomic $\cL_{p}$-formula or the negation of an atomic $\cL_{p}$-formula. \qed
\elem
Then we get the analog of \cite[Lemma 9.2]{BW} with the later remark on relational languages  (pages 305-306).
\lem \label{nota*} For any $\cL$-formula $\varphi(\bar v)$ in prenex normal form: $Q_{1}u_{1}\ldots Q_{m}u_{m}\;\,\psi(\bar u,\bar v)$ where $Q_{i}\in \{\exists, \forall\}$, $1\leq i\leq m$, $\bar u:=(u_{1},\ldots,u_{m})$, $\psi$ is an open $\cL$-formula containing no relation symbols, there is an atomic $\cL_{p}$-formula $\hat{\psi}(w,\bar u,\bar v)$ such that for any element $\cD$ of $\cC$ of cardinality $>1$, for any $\bar b\in \cD$,
\[
\cD^p\models (\forall w\,Q_{1}u_{1}\ldots Q_{m}u_{m}\,\hat{\psi}(w,\bar u,\bar b))\leftrightarrow \cD\models \varphi(\bar b).
\]
We denote $\forall w\,Q_{1}u_{1}\ldots Q_{m}u_{m}\,\hat{\psi}(w,\bar u,\bar v)$ by $\varphi_p(\bar v)$.
\par Assuming now that $\cC$ is the class of models of a complete $\cL$-theory. For any $\cL$-formula $\varphi(\bar w)$ in prenex normal form: $\forall \bar u\exists \bar v\;\,\psi(\bar u,\bar v,\bar w)$ where $\forall$ (respectively $\exists$) means a string of quantifiers $\forall$ (respectively $\exists$), $\psi$ is an open $\cL$-formula (possibly containing relation symbols), there is a conjunction of atomic $\cL_{p}$-formula $\hat{\psi}(z,\bar u,\bar v,\bar w)$ such that for any element $\cD$ of $\cC$ of cardinality $>1$, for any $\bar b\in \cD$,
\[
\cD^p\models (\forall z\,\forall \bar u \exists \bar v\,\hat{\psi}(z,\bar u,\bar v,\bar b))\leftrightarrow \cD\models \varphi(\bar b).
\]
We denote $\forall z\,\forall \bar u \exists \bar v\,\hat{\psi}(z,\bar u,\bar v,\bar w)$ by $\varphi_p(\bar w)$.\qed
\elem

\section{Boolean products of dense pairs}
\par Let $P$ be a new unary relation symbol and suppose that $\cL$ contains at least one constant $c$. Denote by $\cL_{P}:=\cL\cup\{P\}$. Let $\cC$ be a class of $\cL$-structures. Let $\cC_{P}$ be the class of pairs $(\cA,\cD)$ of elements of $\cC$ with $\cD$ an $\cL$-substructure of $\cA$. We view the elements of $\cC_{P}$ as the expansions of elements of $\cC$ in $\cL_{P}$ with the predicate $P$ interpreted by a proper $\cL$-substructure. We will say that the pair $(\cA,\cD)$ is elementary if $\cD\prec \cA$.
\lem \label{pair}Let $\cA\in\Gamma^a(\cC_{P})$ with $\cA=\prod^s_{x\in \cX} \cA_{x}$ and $\cA_{x}\in \cC_{P}$. Define $D:=\{(a_{x})\in \prod^s_{x\in \cX}\cA_{x}\colon \cA_{x}\models P(a_{x})\}$.
Then $\cD\in \Gamma^a(\cC)$ (with $\cX(\cD)=\cX$) and whenever $\cA\in\Gamma^e(\cC_{P})$, $\cD\in\Gamma^e(\cC)$.
\elem
\pr Let $\cD_{x}:=P(\cA_{x})$,  let us first show that $\cD=\prod^s_{x\in \cX} \cD_{x}$. So for $x\in \cX$ and $d_{x}\in D_{x}$, we find $a\in A$ such that $a(x)=d_{x}$ and for all $y\in \cX$, $P(a(y))$. 
By hypothesis, there is $a\in A$ such that $a(x)=d_{x}$. Consider the atomic formula $P(u)$. Then by condition (1) in Definition \ref{bp}, $[P(a)]$ is a clopen set $U\subset \cX$. We have that $[P(a)]\neq \emptyset$. Let $c\in A$ (the interpretation of $c$ in $\cA$) with $c(x)$ the interpretation of $c$ in $\cA_{x}$. Since $\cD_{x}$ is an $\cL$-substructure of $\cA_{x}$, $c(x)\in \cD_{x}$. By the patchwork property (of $\cA$), there is $h\in A$ with $U\subseteq [h=a]$ and $X\setminus U\subseteq [h=c]$. By construction $P(h)$ holds.
\par Now let us show that $\cD\in \Gamma^a(\cC)$. Checking (P1) in Definition \ref{bp} is straightforward. For (P2) in Definition \ref{bp}, let $U$ be a clopen subset of $\cX$, let $f, g\in D$. So there is $h\in A$ such  that $U\subseteq[f=h]$ and $X\setminus U\subseteq[g=h]$.
Since for all $x\in \cX$, $P(h(x))$, we have $h\in D$.
\par Now assume that $\cA\in\Gamma^e(\cC_{P})$. In order to check (P3), consider $\varphi(\bar u)$ an $\cL$-formula and $\bar f$ an $\vert \bar u\vert$-tuple of elements of $\cD$. Let $[\varphi(\bar f)]^{\cD}:=\{x\in \cX: \cD_{x}\models \varphi(\bar f(x))\}$ and denote $\varphi^P$ the formula gotten from $\varphi$ when relativizing the quantifiers to $P$. Note that $\cD_{x}\models \varphi(\bar f(x))$ is equivalent to $\cA_{x}\models \varphi^P(\bar f(x))$. By hypothesis, $[\varphi^P(\bar f)]^{\cA}$ is a clopen subset of $\cX$, and by the above it is equal to $[\varphi(\bar f)]^{\cD}$. 
\qed
\rem \label{el-pair}Let $\cA=\prod^s_{x\in \cX} \cA_{x}$, $\cA_{x}\in \cC$. Let $\cD_{x}\in \cC$ be an $\cL$-substructure of $\cA_{x}$. Let $\cD:=\{(a_{x})\in \cA\colon a_{x}\in \cD_{x}\}$. Suppose that for each $x\in \cX$, $(\cA_{x},\cD_{x})$ is an elementary pair of elements of $\cC$, then for any $\cL$-formula $\varphi(u_{1},\ldots,u_{n})$ and any $n$-tuple $\bar f\in \cD$,  the set $[\varphi(\bar f)]^{\cD}=\{x\in \cX: \cD_{x}\models \varphi(\bar f(x))\}=\{x\in \cX: \cA_{x}\models \varphi(\bar f(x))\}=[\varphi(\bar f)]^{\cA}$.
\par So, if $\cA\in \Gamma^e(\cC)$, then so is $\cD$.
\erem
\bigskip
\par Now assume that $\cC$ is the class of models of a complete $\cL$-theory $T$ (extending the theory of integral domains) with $T$ a geometric $\cL$-theory on the integral domain sort. Recall that a geometric theory is a theory where the model-theoretic algebraic closure $\acl$ satisfies the exchange property (1)  and where the quantifier $\exists^{\infty}$ is eliminated (2). (Note that the first property (1) implies that $T$ eliminates the quantifier $\exists^{\infty}$ (on the domain sort) \cite[Lemma 3.47]{F}.
\par Denote by $\dim_{\acl}$ be the dimension function on definable sets in models of $T$ induced by $\acl$ and we further assume that $\dim_{\acl}$ is a fibered dimension function. 
 \par Let $\cA\models T$, then a definable subset $B\subset A$ is dense in $\cA$  if for any definable subset $Z\subset A$ with $\dim_{\acl}(Z)=1$, then $Z\cap B\neq \emptyset$ \cite[Definition 7.1]{F}. 
\par Let $\cA\in \cC$, then $(\cA, P(A))$ is a dense pair if $P(A)$ is acl-closed and dense in $\cA$. It implies that $P(A)$ is the domain of an elementary substructure of $\cA$  \cite[Lemma 7.4]{F}. Let $T_{P,d}$ be the theory of dense pairs of models of $T$ and let $\cC_{P,d}$ be the class of dense pairs of elements of $\cC$, considered as $\cL_{P}$-structures.
\fct \cite[Theorem 8.3]{F} The theory $T_{P,d}$ is complete.
\efct
\par (In fact in \cite{F}, one places oneself in the more general context of complete theories extending the theory of integral domains with an existential matroid.)
\medskip
\cor \label{cp} Assume that $\cC$ is the class of models of a complete geometric theory $T$ (extending the theory of integral domains) and that $\tilde T$ is a complete theory of Boolean algebras. Let $\cC_{P,d}$ be the class of models of $T_{P,d}$. Then $Th(\Gamma^e_{\tilde T}(\cC_{P,d}))$ is a complete $\cL_{P}$-theory.
\ecor
\pr We apply \cite[Theorem 4.5 (c)]{BW} and the fact that $T_{P,d }$ is complete.\qed 
\lem \label{ep}  Assume that $\cC$ is the class of models of a complete geometric theory $T$ (extending the theory of integral domains) and let $\cC_{P,d}$ be the class of dense pairs of elements of $\cC$. 
Let $\cA\in \Gamma^e(\cC_{P,d})$ and
 consider the expansion $(\cA,P(\cA))$ with $P(A):=\{(a_{x})\colon \cA_{x}\models P(a_{x})\}$. Then  $(\cA,P(\cA))$ is an elementary pair of elements of  $\Gamma^e(\cC)$.
\elem
\pr Let $\cD:=P(\cA)$ and $\cD_{x}:=P(\cA_{x})$. By Lemma \ref{pair}, $\cD\in \Gamma^e(\cC)$ (and $\cX(A)=\cX(D)$).
 Let us check that $\cD\preceq \cA$. Let $\varphi$ be an $\cL$-formula with determining sequence: $(\Phi^*(z_{1},\ldots,z_{\ell}),\psi_{1}(\bar u),\ldots,\psi_{\ell}(\bar u))$ (see Fact \ref{det}). Let $\bar f\in D$. Then $\cA\models \varphi(\bar f) \leftrightarrow \cX^*\models \Phi^*([\psi_{1}(\bar f)],\ldots,[\psi_{\ell}(\bar f)])$. Since for every $x\in \cX$, $\cD_{x} \preceq \cA_{x}$, we have for each $1\leq i\leq \ell$ that $\cA_{x}\models \psi_{i}(\bar f(x))\leftrightarrow \cD_{x}\models\psi_{i}(\bar f(x))$. So $\cA\models \varphi(\bar f) \leftrightarrow \cD\models \varphi(\bar f)$. \qed

\section{Dense pairs of boolean products}\label{top}
\par In this section we specialize to the case when the geometric theory $T$ considered previously is a complete $\cL$-open theory of topological fields of characteristic $0$. In \cite[Definition 1.2.1]{CP}, we dealt with many-sorted structures, but  note that the results on pairs of structures recalled in the previous section, were stated for one-sorted structures, so in the present context, it should be understood to be restricted to the field sort. On the field sort, the language $\cL$ is a relational expansion of the language of fields, namely the language of rings with a multiplicative inverse${\;}^{-1}$ with the convention $0^{-1}=0$, together with a set of constants \cite[Preliminaries 1.1]{CP}. These restrictions on the language imply that the field algebraic closure coincides with the model-theoretic algebraic closure  in models of $T$ \cite[Proposition 1.3.4]{CP}.
\par Models $\cK$ of $T$ are endowed with a definable topology (namely a basis of neighbourhoods of $0$ is given by $\{\chi(K, b)$ with $b$ a tuple varying in $K\}$). 
(On cartesian products of $K$, one puts the product topology).  The theory $T$ is a geometric $\cL$-theory (on the field sort) \cite[Proposition 1.3.4]{CP}
and the topological dimension coincides with the model-theoretic dimension (induced by $acl$) and also with the algebraic dimension coming from the field structure \cite[Proposition 1.3.6]{CP}.

We will consider boolean products of models of $T$ putting condition $(\dagger)$ on relation symbols (see section \ref{projector}). Note that the domain of such structure (on the field sort) is a commutative von Neumann regular ring. 
\par Let us quickly recall some basic facts about these rings. A commutative ring $(R,+,-,\cdot,0,1)$ is von Neumann regular if it satisfies $\forall x\exists y\;(xyx=x)$. 
Let $\cB(R)$ be the Boolean subring consisting of the idempotents of $R$. Let $M_{B}$ be a maximal ideal of $\cB(R)$. Then $M_{B}R$ is a maximal ideal of $R$ and given $M$ a maximal ideal of $R$, we have that $M\cap \cB(R)$ is a maximal ideal of $\cB(R)$ and $M=(M\cap \cB(R))R$.
Let $X$ be the set maximal ideals of $\cB(R)$; when $X$ is viewed as a topological space, we denote it by $\cX(R)$ (as before). Using Stone duality, we have that $\cX(R)^*$ is isomorphic to $\cB(R)$, which is a definable in $R$. 
\par It is easy to see that $R=\prod^s_{x\in X} R/xR$ and furthermore this subdirect product satisfies properties (P1) and (P2) of Definition \ref{bp}. So $R\in \Gamma^a(\bC)$ with $\bC$ the class of fields.
\medskip
\par Now let us consider $\cA=\prod_{x\in X}^s \cA_{x}\in \Gamma^e(T)$ with $\cA_{x}\models T$, as an $\cL_p$-structure (see section \ref{projector}). Note that we have the constant $1$ in the language , so we can express (in $\cA_{x}$) that a term $t$ is different from $0$ by the atomic formula $p(1,t)=0$. Another easy (but useful) remark is that in the class of von Neumann regular rings, the expansion $\cL_{p}$ is an expansion by definitions of $\cL$.
Indeed in any boolean product of integral domains, we can define $p(a,b)$ as follows: 
\begin{equation}\label{eq1}
p(a,b)=c\leftrightarrow \exists d\;(b\,d\,b=b\,\wedge\, b\,c=0\wedge\, (c-a)\,(1-b\,d)=0).
\end{equation}
(One expresses that the supports of $b$ and $c$ are disjoint and on the complement of the support of $b$, $c$ is equal to $a$.) Moreover the defining formula is a (positive primitive) existential $\cL$-formula. Since it defines a function it will imply transfer of model-completeness results from $\cL_{p}$ to $\cL$. 
\par We want first to put a topology on $A$ using the $\cL$-formula $\chi$.  We will assume that the formula $\chi$ is equivalent to a positive primitive existential formula  
$\cL_{p}$-formula $\chi_{p}$ (see Lemma \ref{nota*}), in order to have the following property: for all $a\in A$ and tuple of parameters $b$, $\cA\models \chi_{p}(a,b)\leftrightarrow [\chi(a,b)]=X(\cA)$.
\par We put the following definable topology on $A$: a basis of neighbourhoods of $0$ is given by $\cV:=\{\chi_{p}(A, b)$ with $\chi(0,b)$ and $b$ a tuple varying in $A\}$ and a basis of neighbourhoods of $r\in A$ is of the form $r+V$ with $V\in \cV$. Note that these subsets vary the neighbourhoods of $r$ in the induced topology on the direct product of these topological fields $\cA_{x}$ since for each $x\in X$, $r(x)+\chi(A_{x},b_{x})$ is a neighbourhood of $r(x)\in A_{x}$ by choice of $\chi$.

\par Let us check that indeed this is a Hausdorff topology and that the ring operations are continuous.
\par The ring operations are continuous since it holds in each $\cA_{x}$, the topology being definable, it can be expressed by a formula whose truth value is a clopen subset of $X$ and finally we apply the compacity of $X$ and the patchwork property of $\cA$.
\par Let us show this is Hausdorff. So let $r\neq s\in A$; it suffices to show that there is a neighbourhood of $0$ not containing $r-s$. Let $y\in X$ be such that $(r-s)(y)\neq 0$.
Let $b_{y}$ be a tuple of elements of $A_{y}$ such that $\chi(A_{y},b_{y})$ is a neighbourhood of $0$ (in $A_{y}$) and does not contain $(r-s)(y)$. Let $b$ be a tuple of elements of $A$ such that $b(y)=b_{y}$.
Consider the set of $x\in X$ such that $\chi(0,b(x))\wedge \neg\chi((r-s)(x),b(x))$. This is a clopen subset $U$ of $X$ containing $y$.
Then for each $z\in X\setminus U$, choose a tuple $f_{z}$ of elements of $A$ such that $\chi(A_{z},f_{z}(z))$ is a neighbourhood of $0$. Then we use the fact that $X$ is compact and that $A$ has the patchwork property.
\medskip
\par Now let us consider pairs of models of $\Gamma^e(T)$. We first introduce the following notions. Recall that the language $\cL$ contains at least one constant.
\dfn \label{def:subst} Let $\cA\models \Gamma^e(T)$ and $\cD$ be an $\cL$-substructure of $\cA$ ($\cD\subseteq_{\cL} \cA$). For each $x\in X(\cA)$, let $D_x:=\{u\in A_x: \exists d\in D\;d(x)=u\}$. It is easily checked that $D_x$ is the domain of a substructure of $A_x$ that we denote by $\cD_x$ and that $\cD=\prod_{x\in X(\cA)}^s \cD_{x}$. We will say that $\cD$ is $\acl$-closed in $\cA$ if for every $x\in X(\cA)$, $\cD_x$ is $\acl$-closed in $\cA_x$.
\par Let $\varphi(\bar z)$ be an $\cL$-formula and $\bar f$ be a tuple of elements of $D$. As in Remark \ref{el-pair}, let $[\varphi(\bar f)]^{\cD}=\{x\in X(\cA): \cD_{x}\models \varphi(\bar f(x))\}$. If $\varphi(\bar z)$ is an atomic formula then $[\varphi(\bar f)]^{\cA}=[\varphi(\bar f)]^{\cD}$. 
\par Now we will further assume that 
\begin{enumerate}
\item for any $\cL$-formula $\varphi(\bar z)$ and any tuple $\bar f$ of elements of $D$, $[\varphi(\bar f)]^{\cD}\in X(\cA)^*$,
so the set $\{[\varphi(\bar f)]^{\cD}: \bar f\subset D, \varphi$ an $\cL$-formula $\}$ is a boolean subalgebra of $X(\cA)^*$, that we denote by $X(\cD)^*$, 
\item the patchwork property holds in $\cD$ with respect to $X(\cD)^*$.
\end{enumerate}
If both conditions hold, we will say that $\cD\models  \Gamma^e(T)$ with respect to $X(\cD)^*$.
\edfn
\dfn \label{def:densesubst} Let $\cA\models \Gamma^e(T)$, $\cD\subseteq_{\cL} \cA$ and suppose that $\cD\models  \Gamma^e(T)$ with respect to $X(\cD)^*$. Then the pair $(\cA,\cD)$ is a dense pair if 
\begin{enumerate}[label={(D\arabic*)}]
\item $\cX(\cD)^*\preceq \cX(\cA)^*$,
\item $\cD$ is $\acl$-closed in $\cA$,
\item for every tuple $b$ in $A$, $\chi_{p} (A,b)\cap D\neq \emptyset$. 
\end{enumerate}
\begin{enumerate}[label={(D\arabic*)}, resume]
\item $\forall x\in \cX(\cA)\;\forall e\in \cX(\cA)^*\;\exists \tilde e\in  \cX(\cD)^*\;(e(x)\neq 0\rightarrow \tilde e\subset e\wedge \tilde e(x)\neq 0)$.
\end{enumerate}
\edfn
\par Note that this class of dense pairs is elementary.
\lem \label{elem} Let $\cA\models \Gamma^e(T)$, $\cD\subseteq_{\cL} \cA$ and suppose that $\cD\models  \Gamma^e(T)$ with respect to $\cX(\cD)^*$. Suppose that the pair $(\cA,\cD)$ is a dense pair, then $\cD\preceq \cA$.
\elem
\pr For each $x\in X(\cA)$, we defined $D_x=\{u\in A_x: \exists d\in D\;d(x)=u\}$. By hypothesis on $\chi$, for each $x\in X$, $\chi(A_x,b(x))\cap D_x\neq \emptyset$. By assumption $\cD_x$ is $\acl$-closed so $(\cA_{x},\cD_{x})$ is a dense pair, which implies that $\cD_x\preceq \cA_x$ \cite[Lemma 7.4]{F}. 
\par Let $\bar f\in D$ and assume that 
$\cA\models \varphi(\bar f) \leftrightarrow \cX(\cA)^*\models \Phi^*([\psi_{1}(\bar f)],\ldots,[\psi_{\ell}(\bar f)]),$ where $(\Phi^*, \psi_{1},\ldots,\psi_{\ell})$ is a determining sequence for $\varphi$ (see Fact \ref{det}).
Since for each $x\in X(\cA)$, $\cD_x\preceq \cA_x$, we have that $[\psi_i(\bar f)]^{\cA}=[\psi_i(\bar f)]^{\cD}$, $1\leq i\leq \ell$. So since $X(\cD)^*\preceq X(\cA)^*$, we get that $\cD\preceq \cA$.
\par Note that in the above proof we did not use condition (D4). \qed

\prop Let $(\cA_{0},\cD_{0})\subseteq (\cA,\cD)$ be two dense pairs of models of $\Gamma_{\tilde T}^e(T)$. Suppose the theory $T_{P,d}$ is model complete and  $\neg P(u)$ is equivalent to a positive primitive existential $\cL_{P}$-formula (the condition $(\dagger)$ for $P$). Then $(\cA_{0},\cD_{0})$ is existentially closed in $(\cA,\cD)$.
\eprop
\pr First let us give an equivalent condition for an existential formula to hold in a dense pair of models of $\Gamma_{\tilde T}^e(T)$. Let $\cA\models \Gamma_{\tilde T}^e(T)$. 
Let $\varphi(\bar y)$ be an existential $\cL_P$-formula $\exists \bar u\,\theta(\bar u,\bar y)$, where $\theta(\bar u,\bar y)$ is of the form 
$\bigwedge_{i\in I} \theta_i(\bar u,\bar y)\wedge \bigwedge_{j\in J} \neg\theta_j(\bar u,\bar v)$, where $\theta_{i}, \theta_{j}, i\in I, j\in J$ are atomic $\cL_{P}$-formulas. Note that an atomic formula is either of the form $r(t_{1}(\bar u,\bar v),\ldots,t_{k}(\bar u,\bar v))$, where $r$ is a relation symbol of arity $k$ of $\cL$, or of the form $P(t(\bar u,\bar v))$, or $t(\bar u,\bar v)=0$, where $t(\bar u,\bar v), t_{1}(\bar u,\bar v),\ldots, t_{k}(\bar u,\bar v)$ are $\cL$-terms. For each atomic formula of the form $P(t(\bar u,\bar v))$, we introduce a new variable $w$ and we replace this atomic formula by $P(w) \wedge w=t(\bar u,\bar v)$. So from now on we will assume that the atomic subformulas where there is an instance of the predicate $P$ are of the form $P(w)$, with $w$ a new variable.  Then we decompose the tuple $\bar u$ into two parts: $\bar u_0, \bar u_1$ where a component $u$ of $\bar u$ belongs to $\bar u_0$ if and only if in $\psi_0$ we have an instance of $P(u)$ ($\bar u_0$ can be the empty tuple).
\par Set $\psi_0(\bar u,\bar v):=\bigwedge_{i\in I} \theta_i(\bar u,\bar v)$ and let $\psi_0^+(\bar u,\bar v)$ the formula $\psi_0$ where we removed the atomic subformulas of the form $P(w)$. Consider the set $\cP$ of all non-empty partitions of $J$ and for $(J_1,\ldots,J_\ell)\in \cP$, $1\leq s\leq \ell$, consider the formulas $\psi_{J_{s}}(\bar u, \bar v):=(\psi_{0}(\bar u,\bar v) \wedge \bigwedge_{j\in J_s}\neg\theta_j(\bar u,\bar v))$ and $\psi_{J_{s}}^+(\bar u, \bar v):=(\psi_{0}^+(\bar u,\bar v) \wedge \bigwedge_{j\in J_s}\neg\theta_j(\bar u,\bar v))$. Note that both $\psi_0^+$ and $\psi_{J_s}^+$ are $\cL$-formulas.
\par Let $\varphi_{0}(\bar v):=\exists \bar u \psi_{0}(\bar u,\bar v)$ and $\varphi_{J_s}(\bar v):=\exists \bar u ( \psi_{J_{s}}(\bar u,\bar v))$.

\cl \label{claim} Let $\bar b\in A$, then
\[ (\cA,\cD)\models \varphi(\bar b) \leftrightarrow
\]
\[
\forall x\in X(\cA)\; (\cA_x,\cD_{x})\models \varphi_0(\bar b(x)) \wedge \bigvee_{(J_1,\ldots,J_\ell)\in \cP}\bigwedge_{s=1}^{\ell}\\
 \exists x_s\in  X(\cA) \;(\cA_{x_s},\cD_{x_{s}})\models \varphi_{J_{s}}(\bar b(x_s)).
\]
\ecl
\prcl
$(\rightarrow)$ It is immediate. \\
$(\leftarrow)$
\par First we consider each $x_{s}\in X(\cA)$, $(\cA_{x_{s}},\cD_{x_{s}})\models  \varphi_{J_{s}}(\bar b(x_{s}))$. W.l.o.g. we assume that $(\cA_{x_{s}},\cD_{x_{s}})\not\models  \varphi_{\tilde J}(\bar b(x_{s}))$ for some $J_s\subsetneq \tilde J\subset J$. Let $\bar u^s\in \cA$ with $\bar u^s=\bar u_0^s\bar u_1$ be such that 
$(\cA_{x_{j}},\cD_{x_{j}})\models   \varphi_j(\bar u^s(x_{s}), \bar b(x_{s}))$. Let $\bar d_0^s=(d_1,\ldots,d_m)$ be a tuple of elements of $D$  such that $x_s\in [d_t=u_t]$ for each component $u_t$ of $\bar u_0^s$ (we will denote it by $u_t\in \bar u_0^s$).
The truth value of the $\cL$-formula $[\varphi_{J_s}(\bar u^s, \bar b)]$ is a clopen subset of $\cX(\cA)$ containing $x_s$. 
By hypothesis (on the pair), there is $0\neq \tilde e_{s}\in \cX(\cD)^*$ such that $\tilde e_{s}\subset [\varphi_{J_s}(\bar u^s, \bar b)]\cap \bigcap_{u_t\in \bar u_0} [u_t=d_t]$. 
\par We proceed in the same way for each $x_s$ and get a nonzero idempotent $\tilde e_s\in D$. Then let $\tilde e\in \cX(\cD)^*$ be the union of these idempotents $\tilde e_{s}$.
\par Then we place ourselves on the complement of $\tilde e$, and apply the same procedure for every $x\in X(\cA)$ with $e(x)=0$, considering now the formula $\exists \bar u \varphi_0(\bar u, \bar b(x))$. 
\par So we get a covering of $\cX(\cA)$ by clopen subsets of $\cX(\cD)$ and we use the patchwork property of $\cD$ (respectively $\cA$) with respect to $\cX(\cD)^*$ to get a tuple $\bar d$ (respectively $\bar u_1$) of elements of $\cD$ (respectively $\cA$) such that $\theta(\bar d, \bar u^1,\bar b)$ holds.
\qed
 
\par Then let $\bar b\in \cA_{0}$ and suppose that an existential $\cL_{P}$-formula $\varphi(\bar b$ holds in $(\cA,\cD)$. By the claim above, we have 
\[
\forall x\in X(\cA)\; (\cA_x,\cD_{x})\models \varphi_0(\bar b(x)) \wedge \bigvee_{(J_1,\ldots,J_\ell)\in \cP}\bigwedge_{s=1}^{\ell}\\
 \exists x_s\in  X(\cA) \;(\cA_{x_s},\cD_{x_{s}})\models \varphi_{J_{s}}(\bar b(x_s)).
\]
Now for each $x\in X(\cA)$, we have $(\cA_{0,x},\cD_{0,x})\subset (\cA_x,\cD_{x})$. By assumption, the theory of dense pairs of models of $T$ is model-complete, so we have that $(\cA_{0,x},\cD_{0,x})\preceq (\cA_x,\cD_{x})$. Then it suffices to apply the claim again to get that $(\cA_{0},\cD_{0})\models \varphi(\bar b)$.
\par Note that in the above proof we haven't used condition (D1).
\qed

\section{Generic derivations}
\par In this section we want to consider differential expansion of models of $\Gamma^a(T)$ where again $T$  a complete $\cL$-open theory of topological fields of characteristic $0$ (with the requirements on $\cL$ recalled in the previous section, adding the hypothesis $(\dagger)$ on relation symbols). We expand the models of $T$ with a derivation $\delta$, namely an additive morphism namely: for every $a, b\in R$, $\delta(a+b)=\delta(a)+\delta(b)$, satisfying the Leibnitz rule, namely for every $a, b\in R$, $\delta(ab)=\delta(a)b+a\delta(b)$ and denote the corresponding theory $T_{\delta}$ (in the language $\cL_\delta=\cL\cup\{\delta\}$).
\par Let us denote the theory of the models of $\Gamma^a(T)$ expanded with a derivation $\delta$, $\Gamma^a(T)_{\delta}.$
\par Let $R\models \Gamma^a(T)_{\delta}$ and let $e\in\cB(R)$, then $\delta(e^2)=2\delta(e)e=\delta(e)$. So $\delta(e)(2e-1)=0$ and since $(2e-1)^2=1$, we get that $\delta(e)=0$. 
In particular $\delta$ is trivial on the set $\cB(R)$ of idempotents of $R$.
\par Let $C_{R}:=\{r\in R\colon \delta(r)=0\}$, it is a subring of $R$ and as noted above it contains $\cB(R)$.

\lem Let $(R,\delta)$ be a differential von Neumann regular ring. Let $M$ be a maximal ideal of $R$, then $M$ is a differential maximal ideal of $(R,\delta)$. 
\elem
\pr Let $M$ be a maximal ideal of $R$ and let $r\in M$.
Let us check that $\delta(r)\in M$. Let $s\in R$ be such that $rsr=r$; so $\delta(r)=\delta(r) sr$ since $\delta(sr)=0$. Since $r\in M$, $sr\in M$ and so $\delta(r)\in M$. 
\qed
\cor Let $R$ be a differential von Neumann regular ring. Let $\bC_{\delta}$ be the class of differential fields.
Then $(R,\delta)\in \Gamma^a(\bC_{\delta})$.
\ecor
\pr One notes that given an atomic formula $\theta(x)$ in the language of differential rings, it is equivalent to an atomic formula $\tilde \theta$ in the language of rings in $x, \delta(x),\ldots,\delta^m(x)$, for some $m\geq 0$. \qed
\medskip
\par Now denote by $\bC_{\delta}$ the class of models of $T_{\delta}$.
Similarly we have any differential expansion of a model $\Gamma^a(\bC)$ belongs to $\Gamma^a(\bC_{\delta})$, by \cite[Lemma-Definition 2.2.1]{CP}. 

\medskip
 \par In \cite{CP}, we described a theory $T_\delta^*$ consisting of $T_{\delta}$ together with a scheme of axioms (denoted by (DL)) asserting that if a differential polynomial in one variable of order $m\geq 1$ has an algebraic solution which does not annihilate the separant of that polynomial, then it has a differential solution close to that algebraic solution \cite[Definition 2.2.2]{CP}. 
\par We showed that any model of $T_\delta$ embeds in  a model of $T_\delta^*$ provided the following property $(\dagger)_{\ell arge}$ holds in models $\cK$ of $T$: $\cK$ has an elementary extension $\cK_{0}$ equipped with an henselian valuation such that the valuation topology coincides with the original topology. (Note that it implies that any model of $T$ is a large field since henselian fields are large and being large is an elementary property.)
\par It is straightforward to show that the subfield of constants in a model of $T_{\delta}^*$ is dense (both in the topological sense and according to the model-theoretical definition given above).

\par When $T$ admits quantifier elimination (on the field sort), then $T_\delta^*$ axiomatizes the existentially closed models of $T_\delta$ (on the field sort) \cite[Theorem 2.4.2]{CP} (one shows that $T_\delta^*$ admits quantifier elimination in $\cL_\delta$ on the field sort)).
A consequence of that last result is that the theory $T_{\delta}^*$ is complete \cite[Corollary 2.4.7]{CP} ($T$ was assumed to be complete), and so one can relate the theory of dense pairs and of these differential expansions \cite[Section 4]{CP}.
\par From now on when we consider $T_{\delta}^*$, we will always assume that $T$ satisfies $(\dagger)_{\ell arge}$.
\medskip
\par In \cite{CP}, the following result was stated for one-sorted structures since the expansion $T_{P,d}$ refers to what happens on the field sort.
\fct \cite[Lemma 4.2.1]{CP} \label{constant} Let $\cK_{\delta}\models T_{\delta}^*$, then $(\cK,C_{\cK})\models T_{P,d}$.
 \efct

\medskip
\nota \label{delta} Let $x:=(x_{1},\ldots,x_{n})$, let $\bar m=(m_{1},\ldots,m_{n})\in \bN^n$. Denote by $\bar \delta^{\bar m}(x)$ the tuple $(\bar \delta^{m_{1}}(x_{1}),\ldots, \bar \delta^{m_{n}}(x_{n}))$, with $\bar \delta^{m_{i}}(x_{i})=(x_{i},\delta(x_{i}),\ldots,\delta^{m_{i}}(x_{i}))$, $1\leq i\leq n$.
\enota
\fct \cite[Corollary 2.4.8]{CP} In $T_{\delta}^*$, any $\cL_\delta$-formula $\varphi(x)$ can be put in the form $\psi(\bar \delta^m(x))$ for some $\cL$-formula $\psi$ and $m\geq 0$. 
\efct

\par Assume now that $T$ admits quantifier elimination in $\cL$ (with $\cL$ satisfying $(\dagger)$ (on the field sort). So $T_\delta^*$ admits quantifier elimination in $\cL_\delta$ and in particular has a $\forall\exists$ axiomatisation $\Sigma(T_\delta^*)$ (on the field sort). It enables us to use \cite[Theorem 10.7 (b)]{BW} and the remark on the case when the language contains relation symbols  (pages 305-306 in \cite{BW}). So the class existentially closed  $\cL_{\delta}$-expansions of differential von Neumann regular rings $(R,\delta)$ such that if $x$ is a maximal ideal of $\cB(R)$, the $\cL$-structure $R/xR\models T$, is elementary. Let $T_{at}$ be the theory of atomless boolean algebras.
As a corollary  of the Burris and Werner result and that the projector $p$ can be defined by an existential $\cL$-formula (see equation\ref{eq1}), we obtain:
\cor \label{mc-delta} Assume that $T$ admits quantifier elimination (on the field sort). Then the theory $\Gamma_{T_{at}}^e(T_{\delta}^*)$ is model-complete in $\cL_{\delta}$. \qed
\ecor

In \cite[Theorem 10.7 (b)]{BW}, the theory $\Gamma_{T_{at}}^e(T_{\delta}^*)$ was given an explicit axiomatisation (assuming that the relations in the language satisfy the condition $(\dagger)$). One proceeds as follows.
Given  $\varphi$ a $\forall\exists$ $\cL$-sentence, using the procedure described in Lemma \ref{nota*}, one replaces it  by $\varphi_p$ a  $\forall \exists$ $\cL_p$-sentence. We denote by $\Sigma(T_\delta^*)_p$ the $\cL_p$-sentences obtained from $\Sigma(T_\delta^*)$.
\par Note that for $\cA\in \Gamma^e(T_\delta^*)$, then $\cA\models \varphi$ if and only if $[\varphi_p]=X(\cA)$.
 In our specific setting, we can describe an axiomatisation as follows:
\begin{enumerate}[label=(\Alph*)]
\item the set of axioms expressing that $R$ is a commutative von Neumann regular ring endowed with a derivation $\delta$, with no minimal idempotent,
\item the defining axiom for the projector $p(.,.)$, namely \\$\forall a\forall b\exists d\;\;(b d)b= b\;\wedge\;(p(a,b)-a)(1-b d)=0 \;\wedge\;p(a,b) b=0)$,
\item the set $\Sigma(T_\delta^*)_p$.
\end{enumerate}

\par However we would like here to give a more geometric axiomatization, replacing (C) by $\Sigma(T)_p$ where now $\Sigma(T)$ is a $\forall\exists$ axiomatisation for $T$ and adding the scheme $(G)$ described below.
 First let us begin by recalling some notation and with the following lemma. 
 \par Let $K\{y\}$ be the ring of differential polynomials in one variable. Let $q(y)\in K\{y\}$ of order $m\geq 1$ and $s_q$ be the separant of $q$ (namely the formal derivative of $q$ with respect to $\delta^m(y)$; denote by $q^*(y)$ the ordinary polynomial associated with $q(y)$ in variables $y_0,\ldots,y_m$ and by $s_q^*$ the ordinary polynomial associated with $s_q$.
 \medskip
 \par For $S\subset K^{m+1}$, let $\pi_{m}(S)$ the projection in $K^m$ onto the first $m$-coordinates.
\lem Let $\cK\models T_{\delta}^*$ and let $q(y)\in K\{y\}\setminus\{0\}$ with $\vert y\vert=1$ and $\ord_y(q)=m\geq 1$. Let $S:=\{\bar a\in K^{m+1}\colon q^*(\bar a)=0 \land s_q^*(\bar a)\ne 0\}$. Then $\pi_m(S)$ contains an open set. 
\elem
\pr W.l.o.g. $S\neq \emptyset$ and suppose that for some tuple $\bar a\in K^{m+1}$, $q^*(\bar a)=0 \land s_q^*(\bar a)\ne 0$. Since in $T$ any definable set is a finite union of a Zariski closed set and an open set \cite[]{CP}, if $\pi_m(S)$ does not contain an open set, then $\pi_m(S)$ is a finite union of Zariski closed subset of $K^m$, say $Z$. (We use the following argument which may be found in \cite[Proposition 2.3.2]{CP}.)
\par Let  $\cK_{0}$ be an  $\vert K\vert^+$-saturated elementary extension of $\cK$. Let $v$ be an henselian valuation in $K_{0}$ which coincides the definable topology $\tau$ given in models of $T$. Choose $t_0,\ldots, t_{m-1}\in K_{0}$ with $v(t_{m-1})>>\ldots>>v(t_0)>v(K)$. Let $t=(t_0,\ldots,t_{m-1})$, let $w:K(t)\setminus\{0\}\to \Z^m$ be a coarsening of $v$ which is trivial on $K$ and with $w(t_{m-1})>>\ldots>>w(t_0)>0$. So $w(q^*(a_0+t_0,\ldots,a_{m-1}+t_{m-1},a_m))>0$ and $w(s_q^*(a_0+t_0,\ldots,a_{m-1}+t_{m-1},a_m))=0$. Consider $F:=K(t)^h$ an henselization of $K(t)$ inside $K_0$. So there is $b\in F$ such that $w(a_m-b)>0$ and $q^*(a_0+t_{0},\ldots,a_{m-1}+t_{m-1},b)=0 \land s_q^*(a_0+t_0,\ldots,a_{m-1}+t_{m-1},b)\ne 0$. Since $(a_0+t_0,\ldots,a_{m-1}+t_m-1,b)\in K_{0}$, we get that $K\models \exists d_0\ldots\exists d_m\, q^*(d_0,\ldots,d_{m})=0 \land s_q^*(d_0,\ldots,d_{m})\ne 0\,\wedge\, (d_0,\ldots,d_{m})\neq (a_0,\ldots,a_{m})$. Moreover we may find such $(d_0,\ldots,d_m)$ in a prescribed neigbourhood of $\bar a$ and such that $(d_0,\ldots,d_{m-1})\notin Z$, which is a contradiction. \qed
\medskip
\prop The theory $\Gamma_{T_{at}}^e(T_{\delta}^*)$ can be axiomatized by (A), (B), $\Sigma(T)_p$ together with the following scheme of axioms (G): for any $\cR\models \Gamma_{T_{at}}^e(T)_{\delta}$, for any $\cL$-definable subset $S$ in $R^{n+1}$ such that $\pi_{n}(S)$ contains an open set, there is an element $a$ in $R$ such that $\bar \delta^n(a)\in S$.
\eprop
\pr Let $\cR_{\delta}\models \Gamma_{T_{at}}^e(T)_{\delta}$. In particular $\cR=\prod_{x\in X}^s \cR_{x}$, where $\cR_{x}\models T$ and $\cX$ is an atomless boolean space. Let us show that each $\cR_{x}$ satisfies the scheme (DL) (on the field sort). So it will imply that for each $x\in X$, $(\cR_{x},\delta)$ is a model of $T_\delta^*$. So $\cR_{\delta}\models  \Gamma_{T_{at}}^e(T_{\delta}^*)$.
Let $q(y)\in R_{x}\{y\}\setminus\{0\}$ with $\vert y\vert=1$ and $\ord_y(q)=m\geq 1$. Let $S_{x}:=\{\bar a_{x}\in R_{x}^{m+1}\colon q^*(\bar a_{x})=0 \land s_q^*(\bar a_{x})\ne 0\}$. 
Since $R$ is a subdirect product we may assume that all the coefficients of the polynomial $q(y)$ are of the form $r(x)$ for some $r\in R$ and that $\bar a_{x}$ is of the form $\bar a(x)$ with $\bar a$ a tuple of elements of $R$. Let $U$ be a clopen subset of $X(R)$ included in $[s_q^*(\bar a)\ne 0]$. Let $e\in \cB(R)$ with support $U$. Then multiply every coefficient of $q$ and $s_{q}$ by $e$ and denote the obtained polynomials by $e\, q$ and $e\,s_{q}$ respectively.
Then by the preceding lemma, for every $x\in U$, $\pi_m(S_{x})$ contains an open set. Define $S$ as $\{\bar a\in R^{m+1}\colon e \,q^*(\bar a)=0 \land e\, s_q^*(\bar a)\ne 0\}\cap \bar a+\chi^{m+1}(R,b)$.
Then by the axiomatisation (G) we have a differential tuple $\bar \delta^m(b)\in S$ and since $x\in U$, $\bar \delta^m(b(x))\in S_{x}$.\qed
\prop Assume that $T$ admits quantifier elimination (on the field sort). 
Let $\tilde T$ is a complete theory of boolean algebras. Let $\cA$ be a model of  $\Gamma_{\tilde T}^e(T_{\delta}^*)$. Then 
$(\cA, C_{\cA})$ is a dense pair of model of $\Gamma_{\tilde T}^e(T)$ and so an elementary pair.
\eprop
\pr 
For $x\in X(\cA)$, let $C_{A,x}:=\{u\in A_{x}: \exists r\in C_{A}\;\; r(x)=u\}$. 
Let $\varphi(\bar z)$ be an $\cL$-formula and $\bar f$ be a tuple of elements of $C_{A}$. 
Recall that $\cX(C_{\cA})^*:=\{[\varphi(\bar f)]^{C_{\cA}}: \bar f\subset C_{A}$ with $\varphi$ an $\cL$-formula$\}$.
If $\varphi(\bar z)$ is an atomic formula and $\bar f\subset C_{A}$, then $[\varphi(\bar f)]^{\cA}=[\varphi(\bar f)]^{C_{\cA}}$ and since $T$ admits quantifier elimination and $C_{\cA,x}\models T$ \cite[Lemma 4.2.1]{CP}, this holds for every $\cL$-formula. So $\cX(C_{\cA})^*\subseteq \cX(\cA)^*$. Since any clopen subset of $\cX(A)$ is the support of an element of $A$ and that $\cB(A)\subset C_A$, we have equality and let us denote $\cX(C_{\cA})$ by $\cX$.
\par Also, the patchwork property holds in $C_{\cA}$ with respect to $\cX^*$.
So $C_{\cA}\models  \Gamma^e(T)$ with respect to $\cX^*$.
To check that the pair $(\cA,C_{\cA})$ is a dense pair, it remains to check $(D2)$ and $(D3)$ of Definition \ref{def:densesubst} ($(D4)$ is immediate).
\par For (D2), it holds since it is a condition on the fibers, the field of constants of a differential field is relatively algebraically closed and the field algebraic closure coincides with the model-theoretic algebraic closure in models of $T$.
\par For (D3), 
we use the fact it holds in each fiber since $C_{\cA_{x}}$ is (topologically) dense in $\cA_{x}$, the compacity of $\cX$ and the patchwork property together with the fact that truth values of formulas are clopen subsets. 
\par Finally we apply Lemma \ref{elem}. \qed

\medskip

\lem Let $\cC$ be the class of models of $T$, l et $\tilde T$ is a complete theory of Boolean algebras. Let $\cA\in \Gamma_{\tilde T}^e(\cC_{P})$. Then $(\cA, P(\cA))$ has an $\cL_{P}$-elementary extension $(\cA^*,C_{\cA^*})$ with $\cA^*$ a model of $\Gamma_{\tilde T}^e(T_{\delta}^*)$.
\elem
\pr Let $\cA\in \Gamma_{\tilde T}^e(\cC_{P})$. The theory $Th(\Gamma_{\tilde T}^e(\cC_{P}))$ is a complete $\cL_{P}$-theory (Corollary \ref{cp}).
For $\cK_\delta\models T_{\delta}^*$, $(\cK,C_{\cK})\models T_{P,d}$, interpreting $P$ in $\cK$ by $C_{K}$ (see Fact \ref{constant}). The bounded boolean power $\cK[\cX]^*$, with $\cX^*\models \tilde T$, as an $\cL_P$-structure, belongs to  $\Gamma_{\tilde T}^e(\cC_{P})$.
So $(\cA,P(\cA))\equiv_{\cL_{P}} (\cK[\cX]^*,C_{\cK[\cX]^*})$. Then it suffices to apply Keisler-Shelah's theorem and note that $\cK[\cX]^*$ is a model of 
$Th(\Gamma_{\tilde T}^e(\cC_{\delta}))$.
\qed
\medskip
\par Now let us state a model-completeness result for dense pairs, analogous to Corollary \ref{mc-delta}. In section \ref{app}, we will inspect classical examples of such pairs.
\cor\label{mcpairs} Let $\cC$ be the class of models of $T$ and suppose that the class $\cC_P$ of dense pairs of models of $T$ is model-complete in an expansion $\tilde \cL$ by definitions of $\cL$ by relation symbols satisfying $(\dagger)$. Then the class $\Gamma^e_{T_{at}} (\cC_P)$ is model-complete in $\tilde{\cL}_{P}$. 
\ecor
\pr Again we apply \cite[Theorem 10.7 (b)]{BW}, the remark on the case when the language contains relation symbols  (pages 305-306 in \cite{BW}) and that the fact that the projector $p$ can be defined by an existential $\cL$-formula (see equation (\ref{eq1})).
\qed

\medskip
\par In \cite{CP}, we showed that the $\cL_{\delta}$-theory $T_{\delta}^*$ has $\cL$-open core, namely given a model $\cK$ of $T_{\delta}^*$ and a $\cL_{\delta}$-definable set which is open, it is $\cL$-definable.
\par We proved it by associating with a $\cL_{\delta}$-definable set $Y\subset K^n$, an $\cL$-definable set $Z$ (in some cartesian product of $K$) with the following property \cite[Definition 3.1.2]{CP}:
\begin{enumerate}
\item $Y=\nabla_{\bar m}^{-1}(Z)$, and 
\item $\bar Z=\overline{\nabla_{\bar m}(Y)}$, 
\end{enumerate}
where $\bar m=(m_{1},\ldots,m_{n})\in \bN^n$, $\bar Z$ denotes the topological closure of $Z$ and $\nabla_{\bar m}(Y):=\{(\delta^{\bar m}(a)): a\in Y\}$. We called $(Y, Z, \bar m)$ a linked triple. (See \cite[Proposition 3.1.7, Theorem 3.1.11]{CP}. We showed that one can take $\bar m$ to be the order of $Y$ (\cite[Definition 2.4.10]{CP})).
\medskip

\par Let us assume that $T$ admits quantifier elimination in order to have that the theory $\Gamma_{T_{at}}^e(T_{\delta}^*)$ is model-complete in the language $\cL_{\delta}$ (Corollary \ref{mc-delta}). Let $\cA\models \Gamma_{T_{at}}^e(T_{\delta}^*)$ and let $\varphi(\bar y)$ be an existential $\cL_\delta$-formula. By assumption on the language we may assume that $\varphi(\bar y)$ is of the form $\exists \bar u\,\theta(\bar u,\bar y)$, where $\theta(\bar u,\bar y)$ is a conjunction of atomic and negation of $\cL_{\delta}$-formulas, where the negation of atomic formulas are of the form $t\neq 0$, where $t$ is an $\cL_{\delta}$-term. 
\par Consider an existential $\cL_{\delta}$-formula  of the form
 $\varphi(\bar y):=\exists \bar u (\bigwedge_{i\in I} p_i(\bar u,\bar y)=0\wedge \bigwedge_{j\in J} q_j(\bar u,\bar y)\neq 0 \wedge \bigwedge_k r_k(\bar t_k(\bar u,\bar y))$, where $r_k$ is an $\cL$-relation and $\bar t_k$ a tuple of $\cL_\delta$-terms.
\par Set $\psi_0(\bar u, \bar y):=(\bigwedge_{i\in I} p_i(\bar u,\bar y)=0 \wedge \bigwedge_k r_k)$ and $\varphi_0(\bar y):=\exists \bar u \psi_0(\bar u,\bar y)$.
Let $\psi_{j}(\bar u, \bar y):=\psi_{0}(\bar u, \bar y);\wedge\;q_j(\bar u,\bar v)\neq 0)$ and $\varphi_{j}(\bar y):= \exists \bar u\psi_{j}(\bar u,\bar y)$.
Since the class of models of $\Gamma_{T_{at}}^e(T_{\delta}^*)$ is closed under finite products and since any model of $\Gamma_{T_{at}}^e(T_{\delta}^*)$ is existentially closed, we can apply an observation of S. Burris \cite[Theorem 3]{B}, and get that 
\begin{equation}\label{eq2}
\varphi(\bar y)\leftrightarrow \varphi_0(\bar y) \wedge \bigwedge_{j\in J} \varphi_{j}(\bar y).
\end{equation} 
\par So a determining sequence  for $\varphi$ can be chosen as: $(\Phi^*(z_0, z_1,\ldots, z_{\vert J\vert}),\varphi_0,\varphi_j, j\in J)$, where $\Phi^*(z_0, z_1,\ldots, z_{\vert J\vert}):=(z_0=1\wedge \bigwedge_{j\in J} z_j\neq 0)$ (Fact \ref{det}).

\par Let $\bar m_{0}$ be the order of $\varphi_{0}$ and let $\xi_0$ be a $\cL$-formula such that $(\varphi_{0},\xi_{0},\bar m_0)$ is a linked triple.
Similarly for $j\in J$, let $\bar m_{j}$ be the order of the formula $\varphi_{j}$ and let $\xi_{j}$ be a $\cL$-formula such that $(\varphi_{j},\xi_{j},\bar m_{j})$ is a linked triple. 
\par Let $\bar b\in A$, then $\cA\models \varphi(\bar b)$ iff 
\begin{equation}\label{eq3}
\forall x\in X(\cA)\; \cA_x\models \varphi_0(\bar b(x)) \wedge \bigwedge_{j\in J}\\
 \exists x_j\in X(\cA)\; \cA_{x_j}\models \varphi_{j}(\bar b(x_j)).
\end{equation}
We have that $\varphi_0(A_x)=\nabla_{\bar m_0}^{-1}(\xi_0(A_x))$ and $\overline{(\xi_0(A_x))}=\overline{\nabla_{\bar m_0}(\varphi_0(A_x))}$.
Similarly for $j\in J$, we have that $\varphi_{j}(A_x)=\nabla_{\bar m_j}^{-1}(\xi_{j}(A_x))$ and $\overline{\xi_{j}(A_x)}=\overline{\nabla_{\bar m_{j}}\varphi_{j}(A_x)}$.
\par Since $T_\delta$ is complete and since if $\varphi(A)\neq \emptyset$, there  is $x_j\in X(\cA)$ such that $\varphi_{j}(\cA_{x_j})\neq \emptyset$, we have that for every $x\in X(\cA)$, $\varphi_{j}(\cA_x)\neq \emptyset$.
\par Moreover $T_\delta^*$ eliminates $\exists^{\infty}$ \cite[Theorem A.0.6]{CP} and so either for every $x\in X(\cA)$, $\vert \varphi_{j}(\cA_x)\vert$ is finite and of the same cardinality, or $\varphi_{j}\cA_x)$ is infinite.
\par Note that if $\varphi_{j}(\cA_x)$ is finite then $\varphi_{j}(\nabla_{\bar m_{j}}(A_{x}))$ is equal to $\xi_{j}(\cA_x)$ \cite[Lemma 3.1.9]{CP}.
\par To sum up, we get the following result.
\prop Assume that $T$ admits quantifier elimination and let $\cA\models \Gamma_{T_{at}}^e(T_{\delta}^*)$. Given an $\cL_\delta$-formula $\varphi$ with determining sequence $(z_0=1\wedge \bigwedge_{j\in J} z_j\neq 0)$, one can associate finitely many $\cL$-formulas $\xi_0, \xi_{j}, j\in J$, and finitely many tuples $\bar m_0, \bar m_j, j\in J$ of natural numbers with  the property that for every $\bar f\in A$, 
\begin{equation*}\label{eq4}
\cA\models \varphi(\bar f) \leftrightarrow \big{(}\forall x\in X(\cA)\; \cA_x\models \xi_0(\delta^{\bar m_0}\bar b(x)) \wedge \bigwedge_{j\in J}
 \exists x_j\in X(\cA)\; \cA_{x_j}\models \xi_{j}(\delta^{\bar m_j}(\bar b(x_j)))\big{)}
\end{equation*}
and 
\begin{enumerate}
\item $\overline{\varphi_0(\nabla_{\bar m_0}(A))}=\overline{\xi_0(A)}$,
\item for every $j\in J$, $\overline{\varphi_j(\nabla_{\bar m_j}(A))}=\overline{\xi_j(A)}$.
\end{enumerate} \qed
\eprop

\medskip
\par Let $\cG$ be a collection of sorts of $\cL^{eq}$ and let $\cL^{\cG}$ be the restriction of $\cL^{eq}$ to the field sort together with the new sorts in $\cG$.
\fct \cite[Theorem 3.3.3]{CP} Suppose $T$ admits elimination of imaginaries in $\cL^{\cG}$. Then the theory $T_{\delta}^*$ admits elimination of imaginaries in $\cL_{\delta}^{\cG}$.
\efct
\par Then the natural question is what happens for $\Gamma_{\tilde T}^e(T_{\delta}^*)$? This question was answered in a recent preprint of J. Derakhshan and E. Hrushovski \cite{DH}.  
\par Starting with (complete) theories of boolean algebras, L. Newelski and R. Wencel and then R. Wencel obtained the following results. The theory $T_{at}$ admits weak elimination of imaginaries (a former proof due to J. Truss uses the small index property of atomless boolean algebras and $\aleph_{0}$-categoricity) as well as the theory of boolean algebras with finitely many atoms and the theory of atomic boolean algebras \cite{W}.
\par Then letting $\tilde T$ be respectively the theory of atomic boolean algebras and the theory of atomless boolean algebras, if a theory $T_{0}$ admits elimination of imaginaries, then the theory $\Gamma_{\tilde T}^e(T_{0})$ admits weak elimination of imaginaries \cite[Theorems 3.1, 4.1]{DH}.

\section{Application}\label{app}
In this section we give specific examples to which the results described above apply; these are also examples of open theories $T$ of topological fields \cite[Examples 1.2.5]{CP} and all of them fall in the class of dp-minimal fields. Recall that a  dp-minimal field (with possible extra-structure) is either finite, or algebraically closed or real-closed or has a definable non-trivial henselian valuation \cite[Theorem 1.2]{J23}. Further there were classified in Theorems 1.3 up to 1.6 \cite{J23}. In particular a dp-minimal field which is not strongly minimal is endowed with a definable topology.
\par We consider the classical cases of real-closed fields, real-closed valued fields, algebraically closed fields and p-adically closed fields.
In each case we specify a (one-sorted) language in which these theories admit quantifier elimination and in which the corresponding theory of dense pairs is model-complete (except in the case of real-closed valued fields).
\par Let $\cL$ be the language of rings (with identity) and $\cL_{{\;}^{-1}}$ be the language of fields (we define $0^{-1}=0$). 
\par In case of the order topology, the formula $\chi(x,y)$ can be chosen as follows: $\chi(x,y):= (\vert x\vert \leq \vert y\vert\;\;\&\;y$ is invertible), where $\vert x\vert=x$ if $x\geq 0$ and $\vert x\vert=-x$ if $x<0$. For convenience we will replace the relation symbol $\leq$ by the binary function $\wedge$ (interpreted by the infimum of two elements). (We have $\vert x\vert=x\vee (-x)$.) Let $\cL_{\wedge}:=\cL\cup \{\wedge\}$.
\par In case of the valuation topology, the formula $\chi(x,y)$ will be chosen as follows: $\chi(x,y):= (v(x)\leq v(y) \&\;x$ is invertible). In order to have a one-sorted language, we will replace the valuation $v$ by a binary relation symbol with $x\; \divi \;y$ expressing that $v(x)\leq v(y)$. Let $\cL_{\vert}:=\cL\cup\{\divi\}$.
\begin{enumerate}
\item Let $\RCF$ be the $\cL_{\wedge}$-theory of real-closed fields; by a classical result of Tarski, it admits quantifier elimination. 
\item Now consider the cases of valued fields of characteristic $0$.
\begin{itemize}
    \item let $\ACVF_{0}$ be the $\cL_{\vert}$-theory of  algebraically closed valued fields of characteristic $0$, by a classical result of A. Robinson, it admits quantifier elimination,
    \item let $\RCVF$ be the $\cL_{\vert}\cup\{\wedge\}$-theory of real closed valued fields, by a classical result of G. Cherlin and M. Dickmann, it admits quantifier elimination,
 \item let $\pCF_{d}$ be the theory of  $p$-adically closed fields of $p$-rank $d=ef$ (where $e$ is the ramification index and $f$ the residue degree) and $\cL_{v}$ be the language $\cL_{\vert}$ together with $d$ constants $c_{1},\cdots, c_{d}$ and unary predicates $\{P_{n}: n\geqslant 2\}$, where $P_n(x)$ holds iff $\exists y\;x=y^n$. By classical results of A. Macintyre (when $d=1$) and of A. Prestel and P. Roquette, it admits quantifier elimination  \cite[Theorem 5.6]{PR-84}.
\end{itemize} 
\end{enumerate}
\par Since we consider boolean products of such fields, we have to check that for each relation symbol $r$ in our language, we can express $\neg r$ by a positive existential formula (the condition $(\dagger)$, section \ref{projector}). 
Let us begin with $\divi$. In case of $\pCF_{d}$ that $v(x)<v(y)$ is equivalent to $v(\pi x)\leq v(y)$), where $\pi$ is an element of $K$ with smallest strictly positive valuation. Then for the unary relation $P_{n}$, we use the property that 
if $K$ is p-adically closed valued field and $K^*$ the multiplicative group of the field $K$, then $P_n(K^*)$ is a finite index subgroup in $K^*$.
In the case of $d=1$, we can take cosets representatives in $\bN$ \cite[Lemma 4.2]{B}.
\par In case of $\ACVF_{0}$, we introduce a new relation symbol $x\;\Div\;y$ which expresses that $v(x)<v(y)$.
\medskip
\par In all the above cases the theory $T_{P,d}$ of dense pairs has been shown to be complete \cite[Theorem 8.3]{F} and so the theories $\Gamma_{\tilde T}^e(T_{P,d})$ is complete, where $\tilde T$ is a complete theory of boolean algebras (Corollary \ref{cp}). 
\par Now let us examine in which languages the model-completeness of the theory $T_{P,d}$ has been shown. 
First let us introduce the following $n$-ary relation symbols in a pair of fields $(K,P(K))$, $n\geq 2$.
\par Let $D_{n k}$ is an $n$-ary relation symbol which holds on $x_{1},\ldots,x_{n}$ if $x_{1},\ldots,x_{n}$ satisfy a non-trivial polynomial relation with coefficients in $P(K)$ of degree $\leq k$ and 
let $\ell_{n}(x_{1},\ldots,x_{n})$ iff $x_{1},\ldots,x_{n}$ are linearly independent over $P(K)$.
\par As recalled previously the theory of dense pairs of real-closed fields $K$ is model-complete in $\cL_{\wedge}\cup \{P\}\cup \{D_{n k}: n, k\in \bN_{\geq 1}\}$ 
\cite[Theorem 3.6]{R}. 
For convenience we will use instead $\neg D_{n k}$, in other words define $\tilde D_{n,k}$ by $\neg D_{n k}$ (to get that $\neg \tilde D_{n,k}$ is equivalent to a positive existential formula). Further $\neg P(x)\leftrightarrow \tilde D_{1,1}(1,x)$. The model-completeness result (for dense pairs of real-closed fields) has also be shown in the simpler language $\cL_{\wedge}\cup \{\ell_{n}: n\in \bN_{\geq 1}\}$  
\cite[Proposition 1]{D1} (one defines $P(x)\leftrightarrow \neg \ell_{2}(1,x)$). 
\par When the topological field $K$ is endowed with the valuation topology, we have the following results. 
\par In \cite[Corollary 26]{D}, it is shown that for $T=\ACVF$, the theory $T_{P,d}$ is model-complete in the language $\cL_{\vert}\cup\{\ell_{n}; n\geq 2\}$. (One defines $P(x)\leftrightarrow \neg \ell_{2}(1,x)$, but in our case we keep $P$ since we need this property $(\dagger)$ on the relation symbols. Further $\neg P(x)\leftrightarrow \ell_{2}(1,x)$. 
Note that $\neg(\ell_{n}(x_{1},\ldots,x_{n}))$ iff $\exists z_{1}\ldots \exists z_{n} (\bigvee_{i=1}^n z_{i}\neq 0\;\&\;\sum_{i=1}^n z_{i} x_{i}=0\;\&\;\bigwedge_{i=1}^n P(z_i)$).
\medskip
\par Using the same strategy followed for dense pairs of models of $\ACVF_{0}$, one gets the analogous result for  $\pCF$ (for simplicity we state it for the case of rank $1$). (The additional ingredient is that the subgroup of the $n^{th}$-powers in the multiplicative group of the field is an open subset.)
Consider the language $\cL_{\ell,v}:=\cL_{\vert}\cup\{\ell_{n}; n\geq 2\}\cup\{P_{n}; n\geq 2\}$ to which we add the component functions: $\lambda_{n,i}$, $n\geq 2$ and $1\leq i\leq n$ defined as follows: 
\[
z=\lambda_{n,i}(y,x_{1},\ldots,x_{n}) \leftrightarrow
\]
\[
(\ell_{n}(x_{1},\ldots,x_{n})\;\&\;\neg\ell_{n}(x_{1},\ldots,x_{n},y)\;\&\,\exists z_{1}\ldots\exists z_{n}\;(\bigwedge_{j=1}^n P(z_{j})\;\&\;y=\sum_{j=1}^n x_{j}z_{j}\;\&\;z_{i}=z))\;{\rm or}
\]
\[
(\ell_{n}(x_{1},\ldots,x_{n},y)\;\&\,z=0).
\]
Set  $\cL_{\ell,v,\lambda}:=\cL_{\ell,v}\cup\{\lambda_{n,i}: n\geq 2, 1\leq i\leq n\}$.

Then the theory $T_{P,d}$ admits quantifier elimination in $\cL_{\ell,v,\lambda}$ and it is model-complete in the language $\cL_{\ell,v}$. (Then it can be extended to the case of rank $d$.)
\par A similar result should hold for $\RCVF$ but we haven't checked it.
\bigskip
\par So in each cases when $T\in \{\ACVF_{0}, \RCF, \pCF_{d}\}$, since we described an expansion of $\cL$ in which the corresponding $\tilde \cL$-theory $T_{P,d}$ of dense pairs is model-complete, we get that  the theories $\Gamma^e_{T_{at}} (T_{P,d})$ are model-complete in $\tilde{\cL}_{P}$
(Corollary \ref{mcpairs}). 
\par Finally we give an explicit axiomatisation of the theories $\Gamma^e_{T_{at}}(T)$, when $T\in \{\ACVF_{0}, \RCF, \pCF_{d}\}$.
\par In case $T=\RCF$, first let $T_{f}$ is the $\cL_{\wedge}$-theory of lattice-ordered commutative rings with no nonzero nilpotent elements which are in addition an $f$-ring, namely satisfy the universal axiom $a\wedge b=0\rightarrow a b=0$.
Then let $T_{reg}$ be the theory of von Neumann regular rings with no minimal idempotents which are models of $T_{f}$, where any monic polynomial of odd degree has a root and $\forall x\;(x\wedge 0=0\rightarrow (\exists y\;y^2=x))$. 
\par In case of a valuation topology, we first axiomatize von Neumann regular rings which are subdirect products of valued fields of characteristic $0$ as follows.
\par Let $R$ be a commutative von Neumann regular ring and suppose that $R$ is endowed with a binary relation symbol $\divi$ with the following properties:
for convenience we use $O(x)$ to mean that $1\;\divi\;x$,
\begin{enumerate}
\item $O(1)$,
\item $\forall x\forall y\;\;(O(x) \& O(y)\rightarrow (O(x\pm y)\;\&\;O(x y))$,
\item $\forall x\forall y\exists z\;\;(O(x) \& O(y)\;\&O(z)\rightarrow ((y\,z-x)(zx-y)=0)$,
\item $\forall x\exists y\exists z\;(O(y)\&O(z)\& (y-x)(1-xz)=0$.
\end{enumerate}
Then $R$ is a subdirect product of fields $R/x$, where $x$ is a maximal ideal of $R$, equipped with the binary relation $\divi$ with the property that $O(R/x)$ is a valuation ring. In order to satisfy the condition $(\dagger)$, we use the binary relation symbol $x\,\Div\,y$ (interpreted by $v(x)<v(y)$). Set $M(x)$ to mean that $1\,\Div\, x$.
Then we express that $\forall x \forall y\,(O(x)\& M(y)\rightarrow M(x\,y)$ and $\forall x\exists y\;(O(x)\&O(y)\&\neg M(x)\rightarrow M(1-x\,y))$.
\par To express that for every maximal ideal $x$, $R/x$ is a field of characteristic $0$, we say for every $n\in \bN_{\geq 1}$:
$\forall e\;(e^2=e\rightarrow n\; e\neq 0.$
\par Let $T_{v}$ be the theory of von Neumann regular rings endowed with two binary relations symbol $\divi,\;\Div$ satisfying the above. Then let $T_{reg,v,0}$ be the theory of von Neumann regular rings with no minimal idempotents which are models of $T_{v}$ and where any monic polynomial has a root. This corresponds to $T=\ACVF_{0}$.
\par Let $T_{reg,v,p}$ be the theory of von Neumann regular rings with no minimal idempotents which are models of $T_{v}$ and where 
\begin{enumerate}
\item every polynomial of the form $x^n+x^{n-1}+u_{n-2} x^{n-2}+\ldots+u_{0}$, where $M(a_{n-2}),\ldots, M(a_{0})$, has a zero, $n\in \bN_{\geq 1}$ (the henselian property),
\item $\forall x \forall y\;(x\,\Div\, y\rightarrow p\,x\;\divi\; y)$ (the value group is discrete),
\item for every $m\in \bN_{\geq 2}$, $\forall a\exists b\exists e\;\exists e'\;(a\neq 0\rightarrow 
O(e)\wedge O(e')\wedge e e'=1\wedge \bigvee_{\ell=0}^{m-1} a e p^{\ell}= b^m)$ (the value group is a $\Z$-group),
\item $\forall x\,( O(x)\rightarrow M(x (x-1)\ldots(x-(p-1)))$ (the residue field is isomorphic to $\bF_{p}$)
\end{enumerate}
This corresponds to $T=\pCF$ in the rank $1$ case. (Otherwise one replaces $p$ by $\pi$ and expresses that the dimension of $O/p$ over $\bF_p$ is $d$, using the constants $c_1,\ldots, c_d$.)

\end{document}